\documentclass{article}

\usepackage{blindtext}
\usepackage[T1]{fontenc}
\usepackage[utf8]{inputenc} 
\usepackage[british]{babel} 
\usepackage{url}
\usepackage{amsthm, amssymb}         
\usepackage{blindtext}
\usepackage{amsmath}

\theoremstyle{plain}
\newtheorem{theorem}{Theorem}[section]
\newtheorem{corollary}[theorem]{Corollary}

\theoremstyle{definition}

\usepackage{xcolor}

\newcommand{\xa}{\left<x,\alpha\right>}
\newcommand{\ya}{\left<y,\alpha\right>}
\newcommand{\SaR}{\sum_{\alpha\in R_+}}

\newcommand{\cond}[1]{\, #1 \vert \,}
\usepackage{authblk}
\begin{document}

\title{Recent developments in stochastic systems with singular dynamics \\
-- a brief survey --}

\author{
Sergio Andraus\thanks{\texttt{sergio.andraus@japan-iu.ac.jp}},
Nicole Hufnagel\thanks{\texttt{nicole.hufnagel@hhu.de}},
Jacek Małecki\thanks{\texttt{jacek.malecki@pwr.wroc.pl}},
Vlad Margarint\thanks{\texttt{vmargari@charlotte.edu}}
}
\maketitle

\begin{abstract}

Diffusions with singular drifts appear naturally in the study of certain classes of interesting particle systems. Due to the singularities involved in the equations describing them, studying their dynamics is both interesting and challenging as they fall outside of the usual regularity assumptions on the vector fields in the literature. One needs to introduce and study them with different methods, with the use of symmetric polynomials being a particularly promising approach. In this expository work, we give an overview of some models and methods, and we further explore their interactions with the areas of Multiple Loewner chains (in particular, multiple Schramm-Loewner Evolutions), the new links with Random Matrices, and the novel area of Machine Learning for Singular Dynamics. 

\end{abstract}

\section{Introduction}
This paper studies stochastic diffusions with singular drift terms, which fall outside the classical framework where the drift and volatility coefficients are assumed to be at least locally Lipschitz. In such settings, standard existence and uniqueness results for stochastic differential equations do not directly apply, and additional analytical tools are required to understand the behavior of the associated processes. We investigate what structural properties can still be analyzed when these regularity assumptions fail. 

Our approach uses techniques based on symmetric polynomials to study the structure of these diffusions and to obtain insight into their dynamics despite the presence of singularities. In particular, we focus on models in which the drift contains singular interaction terms that prevent particle collisions when the interaction strength parameter is sufficiently large. 

As a central example, we present the Dyson Brownian motion, a system of interacting diffusions with repulsive singular drift that arises in random matrix theory. We discuss how the polynomial framework naturally captures the symmetric structure of the system and allows one to analyze important properties such as the evolution of eigenvalues and the role of singular interactions. We further extend these ideas to study the multiple Schramm-Loewner evolutions model in which the Dyson Brownian motion appears as a driver. 

\section{Interacting particle systems - general case}
Stochastic particle systems have been an important topic in both mathematics and physics for many years. In particular, systems such as the Dyson Brownian motion (DBM) arising in nuclear physics \cite{Dyson62A} and Wishart processes appearing in statistics \cite{Bru91} are well-known examples, due to their deep connections to random matrix theory and integrable systems \cite{KatoriTanemura04, Forrester10}. These particle systems are special cases of the more general particle system $(x(t))_{t\geq0}$ of dimension $N$ given via the stochastic differential equation (SDE)
\begin{equation}
	\label{eq:x:SDE:mult}
	dx = \sigma(x)dB+b(x)dt+\SaR k_{\alpha}(x)\dfrac{\alpha}{\xa} \,dt\/,
\end{equation}
where we slightly abuse the notation by writing $\sigma(x)=(\sigma(x_1),$ $\ldots, \sigma(x_N))$ and $b(x) = (b(x_1),\ldots,b(x_N))$. We denote the Euclidean inner product by $\xa$, and $R_+$ denotes the positive subsystem of the root system $R$, which could be $A_{N-1}$, $B_N$, or $D_N$, see \cite{Humphreys1997} for details. 
The coefficient $k_\alpha(x)$ is the previously mentioned interaction strength parameter, and it is proportional to the $\beta$ parameter from random matrix theory when \eqref{eq:x:SDE:mult} is specialized to the DBM. For simplicity, we suppress the explicit dependence on the time variable $t$. In this paper we concentrate on the main results and therefore do not state all assumptions explicitly. Precise definitions of the objects appearing in the SDE and the corresponding assumptions can be found in \cite{Malecki2026,AndrausHufnagelMalecki25}. 

The first result, which enables us to study more refined properties of the set of collision times, characterizes the nature of collisions, where a collision corresponds to the event $\xa=0$. In such situations, the singular drift term may blow up, and the result ensures that its behaviour remains under control.

\begin{theorem}
    \label{th:multiplecollisions}
    Let $x=(x_1,\ldots,x_N)$ be a solution to \eqref{eq:x:SDE:mult}, assume that $b, k_\alpha, \sigma$ are continuous in their domains, that $\sigma$ and $k_\alpha$ are strictly positive for every root $\alpha$, and only in the $D_N$ case, assume that $\sigma(-u)=\sigma(u)$. Then, no distinct root hyperplanes are hit simultaneously at positive times before the lifetime of the process.
\end{theorem}

In essence, this theorem states that under mild conditions on the coefficients $b, \sigma,$ and $k_\alpha$, all hits are simple, namely, they only involve one root hyperplane at a time. For a deeper understanding of the hitting times, this theorem plays a crucial role. We are interested in the frequency of such hitting times and, in particular, in a quantitative way of measuring it. For this purpose we consider the Hausdorff dimension, a notion of fractional dimension that extends the classical integer-valued notion of dimension to positive real values, thereby allowing a finer distinction between sets, see \cite[8.1 Hausdorff dimension]{Adler1981}. We denote this quantity by $\dim(\cdot)$.

In \cite[Theorem 2, Theorem 3]{AndrausHufnagelMalecki25} we establish, under mild technical conditions on the coefficients, non-trivial upper and lower bounds for $\dim(x^{-1}(\partial W))$, where $\partial W$ is the boundary of the closed Weyl chamber, $\overline W:=\{y\in \mathbb{R}^N \cond{} \exists \alpha \in$ $R_+ : \ya \ge 0\}$. This quantity is naturally a number between $0$ and $1$, since it measures the Hausdorff dimension of a subset of the time axis. It describes the dimension of the set of times at which collisions occur.

\section{Particular cases}
We now discuss several cases covered by our framework.

\subsection{Dyson Brownian motion}

\subsubsection{Hausdorff dimension}
The Hausdorff dimension for the well-known DBM was first calculated in \cite{HufnagelAndraus24}. There, the authors derive bounds for the Hausdorff measure using straightforward and rigorous analytical arguments. In particular, they estimate the density and apply a capacity argument. In \cite{AndrausHufnagelMalecki25}, this result is generalized to processes solving \eqref{eq:x:SDE:mult} by employing symmetric polynomials and using the well-known Hausdorff dimension for the hitting time of zero of a classical one-dimensional Bessel process.

The DBM \cite{Dyson62A} is obtained by setting $\sigma(x)=1$, $b(x)=0$, $k_\alpha(x)=k>0$, and $R_+:=\{e_j-e_i\cond{} 1\leq i<j\leq N\}$. The process is given by the following SDE,
        \[dx_i=dB_i+k\sum_{j=1:j\neq i}^N\frac{dt}{x_i-x_j}.\]
        We note here that $k$ is related to the parameter $\beta$ from random matrix theory by $k=\beta/2$. There are versions of this dynamics where the interaction term is scaled by a factor of $1/N$ \cite{CampbellLuhMargarint25, JiamVlad}. By Corollary~\ref{cor:multiBessel}, we recover the following familiar result \cite{HufnagelAndraus24}: the set of times where particle collisions occur in the DBM has a Hausdorff dimension given by
        \[\dim x^{-1}(\partial W)=\max\Big\{0,\frac12-k\Big\}.\]
    An important feature of this process is that it can be expressed as the eigenvalue process of a real symmetric ($k=1/2$) or complex Hermitian ($k=1$) random matrix with independent Brownian motions up to symmetry \cite{KatoriTanemura04}. These matrix-valued processes are known for showing no collisions, which can be observed in the Hausdorff dimension formula above, as they correspond to $k=1/2$ and $k=1$, that is, the Hausdorff dimension is zero in those cases because there are no collisions, but the converse is not generally true.
    \subsubsection{Applications to Multiple Schramm-Loewner Evolutions}
    This subsection explores the connections to an important development in Mathematical Physics in the last twenty years: the study of Schramm-Loewner Evolutions (SLE). SLE is a family of fractal curves parameterized by $\kappa \geq 0$ introduced by Schramm \cite{Schramm2000} which appears in the context of studying scaling limits of many planar statistical physics models.

One important element in SLE theory is the Loewner Differential Equation. We recall that given a simply-connected domain $D \subset \mathbb{C}$, and a simple curve, $\gamma(0, \infty)$, then for any $t \geq 0$ there exists a conformal isomorphism from $D \setminus \gamma(0, t] \to D$ by the Riemann Mapping Theorem. Focusing on the upper half-plane case. the normalized maps $g_t$ satisfy the Loewner Differential Equation \cite{BerestyckiNorris}: 
\begin{equation}
    \partial_t g_t(z) = \frac{2}{g_t(z) - \sqrt{\kappa}B_t}.
\end{equation}

  One generalization of the SLE model which received attention in the last years is a multi-driver version of the Loewner Differential Equation, which is the Multiple SLE in which the Brownian Motion driver is replaced by a DBM. For works in this direction in the SLE setting, we refer the reader to \cite{HealeyLawler,JiamVlad,EveViv,CampbellLuhMargarint25,Toma}. In addition, the works in \cite{Zhan,Peltola,Dub} cover other models.
This Multiple SLE model has connections with other areas such as Conformal Field Theory (see for example \cite{Cardy}).

The Multiple SLE with Dyson drivers is then defined by setting $k=4/\kappa$ and introducing $\lambda_i:=\sqrt{\kappa}x_i$ to then write the Loewner evolution equation as
\begin{equation}
    \partial_{t} g_{t}(z) = \sum_{i=1}^{N}\frac{2}{g_{t}(z)-\lambda_i(t)}, \quad g_0(z) = z.
\end{equation}
This bridges our notation, $k$, with the canonical SLE theory parameter $\kappa=8/\beta=4/k$.

In our recent work \cite{AndrausHufnagelMaleckiMargarint2026}, we develop these ideas further by introducing a systematic framework for studying multiple Loewner chains driven by interacting particle systems. The key ingredient is a family of symmetric polynomials in the squared distances between the Loewner image of an initial point and the driving particles. Their stochastic dynamics allow us to reduce the boundary-hitting problem to the study of a one-dimensional squared Bessel process. In particular, for the classical Dyson-type models, this yields a precise description of two distinct phase transitions: the driving particles collide when \(\kappa>8\), while, in the absence of particle collisions, the Loewner image hits a driving particle in finite time precisely when \(4<\kappa\leq8\), and when $0<\kappa\leq 4$ there are no hitting events between the Loewner image and the drivers, nor between the drivers themselves.
We further extend this approach to root-system-adapted multiple Loewner chains, and investigate which driver is hit when  there are only two drivers performing a DBM, and the initial point lies between them. These results demonstrate that symmetric polynomials provide a unified framework for understanding boundary behavior in a broad class of multiple Loewner evolutions. 

We currently explore the connections between these techniques and the toolbox that connects SLE and Random Matrix Theory (RMT) introduced in \cite{CampbellLuhMargarint25}. One feature of this direction of research consists of studying the Multiple SLE via a very important tool in RMT, the resolvent
of random matrices. This gives the possibility of an alternative study of Multiple SLE by seeing the right-hand side  of the Loewner equation as the resolvent of a time-evolving random matrix. The evolution of these matrices in time is more regular than the singular evolution of their eigenvalues. Indeed, the dynamics in the matrix space can be of the form 

$$dA_t=-\frac{1}{2}A_tdt+\frac{1}{\sqrt{N}}dB_t.$$ 
Thus, one can consider the corresponding shifted matrix dynamics, compute their resolvents and use them to study the multiple SLE model.

\subsection{Multivariate Bessel process}
Also known as radial Dunkl processes \cite{RoslerVoit98,GallardoYor05}, these are processes given by the SDE
    \[dx_i=dB_i+\SaR k_\alpha\frac{\alpha_i}{\xa}\,dt,\]
    that is, $\sigma=1$, $b=0$, and $k_\alpha(x)=k_\alpha$, with the added constraint that $k_\alpha$ must be invariant under reflections along the roots in $R_+$. A smaller value of $k_\alpha$ increases the probability that the corresponding boundary is hit. In particular, it is well known that collisions can occur only if 
\[
\min_{\alpha \in R_+} k_\alpha < \tfrac{1}{2}, 
\]
see \cite{Demni09}. Assuming that $R$ is one of the classical root systems ($A_{N-1}$, $B_N$, or $D_N$) we have
    \begin{corollary}
        The collision times of the Bessel process have the following Hausdorff dimension,
        \[\dim x^{-1}(\partial W)=\max\Big\{0,\frac12-\min_{\alpha\in R_+}k_\alpha\Big\}.\]
        \label{cor:multiBessel}
    \end{corollary}
    
Hence, the smallest value of $k_\alpha$ dominates and determines the dimension. In particular the case $R=A_{N-1}$ is the DBM.
\subsection{Multivariate Bessel process of type $B_N$}
Of particular interest is the case where $R=B_N$, in which $k_\alpha$ is reduced to two independent parameters, $k_1$ 
    and $k_2$, and $R_+:=\{e_i\cond{}1\leq i\leq N\}\cup\{e_j\pm e_i\cond{}1\leq i<j\leq N\}$. Its SDE reads
    $$dx_i=dB_i+\frac{k_1}{x_i}\,dt+k_2\sum_{j=1:j\neq i}^N \Big\{\frac{1}{x_i-x_j}+\frac{1}{x_i+x_j}\Big\}\,dt.$$
    By Corollary~\ref{cor:multiBessel}, the corresponding Hausdorff dimension of collision times is
    \[\dim x^{-1}(\partial W_{B_N})=\max\Big\{0,\frac12-\min\{k_1,k_2\}\Big\}.\]
    The parameters $k_1$ and $k_2$ control the type of collision that occurs. Indeed, by \cite[Prop.~1]{Demni09}, the first hitting time of zero is finite if $k_1<1/2$, and infinite otherwise, while the first collision time between particles is finite if $k_2<1/2$ and infinite otherwise.
\subsection{Wishart processes}
The Wishart process \cite{Bru91,KonigOConnell01}  is another important system that can be formulated as the matrix eigenvalue process of the product of a rectangular matrix with independent Brownian motion entries with its transpose. As in the Dyson case, the values $k=1$ and $2$ correspond to real and complex matrix formulations of the system, but here we consider the general case given by the SDE
    \begin{align*}
        dy_i&=2\sqrt{y_i}\,dB_i+k a\,dt+k\sum_{j=1:j\neq i}^N\frac{y_i+y_j}{y_i-y_j}\,dt
    \end{align*}
    for $k,a>0$.
    The Wishart process can be obtained from the multivariate Bessel process of type $B_N$ by the relationships $y_i=x_i^2$, $k a=2k_1+2k_2(N-1)+1$, and $k=2k_2$.
    \begin{corollary}
        When $a\geq 2/k +N-1$, the Hausdorff dimension of collision times between particles in a Wishart process is given by
        \[\dim y^{-1}(\partial W_{B_N})=\max\Big\{0,\frac{1-k}{2}\Big\}.\]
    \end{corollary}
    
    \subsection{Jacobi processes}
    Jacobi processes are mutually repelling particle systems in a finite segment of the real line. As introduced in \cite{Demni10}, all particles are ordered and lie in the interval $(0,1)$, but for the sake of symmetry we consider the process on $(-1,1)$, namely, the $N$ particles $\{\lambda_i\}_{i=1}^N$ in \cite{Demni10} are replaced by $x_i=2\lambda_i-1$, and the corresponding SDE reads
    \begin{align*}
        dx_i=\sqrt{1-x_i^2}\,dB_i&+k\bigg[\frac{1}{2}(p-q-(p+q)x_i)\\
        &+\sum_{j=1:j\neq i}^N\frac{1-x_ix_j}{x_i-x_j}\bigg]\,dt,
    \end{align*}
    where $p$ and $q$ are integer parameters satisfying $\min\{p,q\}\geq N-1+2/k$ in order to ensure the SDE has a unique strong solution.
    We note here that the cases $k=1$ and $2$ can be formulated as dynamical extensions of the Jacobi ensembles, very much in the same way as the Wishart and Dyson cases above.
    From the SDE, it is clear that $k>0$ governs the repulsion between particles, while $p$ and $q$ represent the repulsion strength of the left and right walls respectively.

    In this case, we only get a trivial upper bound 
    $$\dim x^{-1}(\partial W)\leq\max\Big\{0,\frac12-k\Big\}.$$

\subsection{Multiple Loewner chain driven by particle systems}

In addition to the aforementioned directions, following the recent developments in \cite{ShrotriMargarint2026} in the one-curve SLE case, it is interesting to consider the analysis of Machine Learning models in the context of Multiple Loewner Chains driven by singular interacting particle systems. The particular case of Dyson Brownian Motion interacting particle systems drivers that leads to multiple SLE, has further interesting links with Theoretical Physics. In \cite{ShrotriMargarint2026} the authors use, between others, Deep Learning architectures Neural Networks to learn the parameter $\kappa \in [0, 8]$ in the one-curve SLE dynamics. The singularity in the dynamics makes the analysis of these models both challenging and interesting. The case of Multiple Loewner Chains adds another layer of complexity as now the drivers themselves have singularities as well. It is thus interesting to study the learning of the parameters in these new extended cases.

\end{document}